\documentclass[12pt]{article}
\usepackage{amsmath, amsfonts, amssymb}
\usepackage{amsmath}
\usepackage{graphicx}
\usepackage{color}
\title{Iterative roots of exclusive multifunctions
\thanks{Supported by NSFC\#11501471(L. Liu), Zhejiang Provincial Natural Science Foundation of China under Grant \#LY18A010017 (L. Li), NSFC\#11821001 (W. Zhang) and NSFC\#11831012 (W. Zhang) } }
\author{
{\sc Liu Liu}\,$^{a}$,~~~ {\sc Lin Li}\,$^{b}$
\footnote{Corresponding author: matlinl@mail.zjxu.edu.cn} ,~~~ {\sc Weinian
Zhang}\,$^{c}$
\\
$^{a}${\small Department of Mathematics, Southwest Jiaotong University}\\
{\small Chengdu, Sichuan 610031, P. R. China}\\
$^b${\small College of Mathematics, Physics and Information Engineering}\\
{\small Jiaxing University, Jiaxing, Zhejiang 314001, P. R. China}
\\
$^{c}${\small School of Mathematics, Sichuan University}\\
{\small Chengdu, Sichuan 610064, P. R. China}
}
\date{}
\begin{document}
\maketitle

\parindent 20pt

\begin{abstract}
In this paper we investigate iterative roots of strictly monotone
upper semi-continuous multifunctions having finitely many jumps. Known results are
concerning roots of order 2 for multifunctions of exact one jump.
For the general investigation, we introduce a concept `intensity' to
formulate the growth of jumps under iteration and find a class of
strictly monotone and upper semi-continuous multifunctions of intensity 1, called
exclusive multifunctions, each of which has an absorbing interval.
Then we use the absorbing interval to construct iterative roots of
order $n$ for those exclusive multifunctions.

\vskip 0.1cm
{\bf Keywords}:
iterative root; multifunction; upper semi-continuity; jump; intensity.

\vskip 0.1cm
{\bf AMS(2000) Subject Classification}: 39B12; 37E05; 54C60

\end{abstract}

\renewcommand{\theequation}{\thesection.\arabic{equation}}
\newtheorem{df}{Definition}[section]
\newtheorem{lm}{Lemma}[section]
\newtheorem{thm}{Theorem}[section]
\newtheorem{pro}{Proposition}[section]
\newtheorem{exmp}{Example}[section]
\newtheorem{rmk}{Remark}[section]
\newtheorem{alg}{Algorithm}[section]
\newtheorem{cor}{Corollary}[section]
\newtheorem{ap}{Application}[section]

\setlength\arraycolsep{2pt}

\baselineskip 16pt
\parskip 10pt


\section{Introduction}
\setcounter{equation}{0}

For an integer $n\ge 1$,
the $n$-th order iterate $g^n$ of a self-mapping $g: X\to X$, where $X$ is a nonempty set, is defined recursively by
$g^n(x)=g(g^{n-1}(x))$ and $g^0(x)=x$.
The iterative root problem for a given self-mapping $f: X\to X$ is to find a self-mapping $g: X\to X$
such that the functional equation
\begin{eqnarray}
g^n=f
\label{itroot}
\end{eqnarray}
is true on $X$. We call $g$ an $n$-th order {\it iterative root} or
{\it fractional iterate} of $f$. Although its research can be
pursued up to the 19th century
(see Babbage's work \cite{Babbage1815}),
this problem has been attractive because it is a weak version of the
embedding flow problem (\cite{Fort,Irwin})
of dynamical systems
and applicable to informatics (\cite{PGLKPM,Kindermann98})
and
was answered systematically for monotone and continuous interval-mappings
(\cite{
Kuczma68, KCG1990, Targonski1981})
around 1960's
because of difficulties caused by the iteration operation.
After 1980's, some results
(see e.g. \cite{Li-Yang-Zhang08,
LiuJarczykLiZhang2011NA,
Liu-Zhang11,
ZhangYang1983AMS}) on construction of
continuous iterative roots were made for interval-mappings having
finitely many non-monotonic points, which destroy the orientation of
iterates.
For cases of dimension $\ge 2$, no results on
iterative roots are found except those obtained in \cite {Lesniak1,Lesniak2}
for planar Sperner homeomorphisms and Brouwer homeomorphisms.

Encountering difficulties in finding iterative roots,
ones also made efforts (see e.g. \cite{Jarc-Powi,Jar-Zh07Ele,LiZhangResult11,
Asmajdor85})
to find multivalued iterative roots,
which have at least one set-valued point (or called jump simply). In general, the composition $G\circ F$ of
multifunctions $F:X\to 2^Y$ and $G:Y\to 2^Z$ is defined by
$
(G\circ F)(x)=G(F(x)),
$
where the image $F(A)$ of a set $A\subset X$ is defined by $F(A):=\cup_{x\in A}F(x)$.
Then, as shown in \cite{Jar-Zh07Ele}, the $n$-th order iterate $F^n$ is defined recursively by
$$
F^n(x)=\bigcup_{x\in X}\{F(y):y\in F^{n-1}(x)\}
\ \ \ \mbox{and}\ \ \ F^0(x)=\{x\}
$$
for all $x\in X$.
One can check that
$$
F^n(A)=\{x_0\in X: \bigvee_{x_1,...,x_{n-1}\in X}\bigvee_{x_n\in A}x_{i-1}\in F(x_i),i=1,...,n\}
$$
for every $n\in \mathbb{N}$ and therefore $f^n(f^m(A))=f^{n+m}(A)=f^m(f^n(A))$.
With these concepts,
Powierza and Jarczyk (\cite{Jarc-Powi})
discussed existence of smallest set-valued iterative roots $g$ of
bijections $f$ in the sense that
$
f(x)\in g^n(x)
$
for all $x\in X$. Another idea (\cite{Jar-Zh07Ele, Li-Jar-Zhang09Publ}) is to find
solutions $g: X\rightarrow 2^{X}$ of (\ref{itroot}) for a
multifunction $f$. We simply refer the former's to the {\it
inclusion sense} but the latter one to the {\it identity sense}. In
\cite{Jar-Zh07Ele} two results on the non-existence of square iterative roots of multifunctions
were obtained on the general nonempty set $X$. References
\cite{Li-Jar-Zhang09Publ,GLY18} are contributed to square iterative roots
of multifunctions with only one jump. For general nonempty set $X$,
the nonexistence of iterative roots of order $2$ is discussed.
For $X=I:=[a,b]\subset \mathbb{R}$, by the method of piecewise construction,
the existence of iterative roots of order $2$
is given for usc (abbreviation of upper semi-continuous) multifunctions.

In this paper we investigate the $n$-th order iterative roots on
$I=[a,b]$ in the identity sense for strictly monotone usc
multifunctions having finitely many jumps. In section 2 we introduce
the concept {\it intensity} to formulate the growth of jumps under
iteration and give properties of strictly monotone usc
multifunctions having finitely many jumps under iteration. We prove
that those multifunctions of
intensity 1,
called {\it exclusive multifunctions},
do not increase the number of jumps under iteration
but have an {\it absorbing interval} each.
In sections 3 and 4,
we use the absorbing interval to construct increasing iterative roots and decreasing ones respectively for increasing exclusive multifunctions.
In section 5, we find iterative roots for decreasing exclusive multifunctions.
We demonstrate our theorems with examples in section 6 and remarks for unsolved cases.


\section{Intensity and absorbing interval}

Let $I:=[a,b]\subset \mathbb{R}$. As defined in \cite{SHNSP}, an point $x_0\in I$ is referred to as
a {\it set-valued point} or simply a {\it jump} of $F$ if $\#F(x_0)\geq 2$.
$F: I\rightarrow 2^I$ is {\it upper
semi-continuous} (abbr. {\it usc}) at $x_0\in I$ if for every open
set $V\subset \mathbb{R}$ satisfying $F(x_0)\subset V$ there exists
a neighborhood $U_{x_0}$ of $x_0$ such that $F(x)\subset V$ for
every $x\in U_{x_0}$. Moreover, we say that $F$ is {\it strictly
increasing} (resp. {\it strictly decreasing}) if $\max F(x)< \min
F(y)$ (resp. $\min F(x)>\max F(y)$) whenever $x,y\in I$ and $x<y$.
Let ${\cal F}(I)$ denote the class of all strictly monotone usc multifunction $F: I\to 2^I$.
For preliminaries, in this section we give properties for iteration
of multifunctions in ${\cal F}(I)$.

\begin{lm}
Let $n\in \mathbb{N}$. If $F\in {\cal F}(I)$ then $F^{n}\in {\cal
F}(I)$. Moreover, if $F$ is increasing then $F^n$ is increasing, and
if $F$ is decreasing then $F^{2n-1}$ is decreasing and $F^{2n}$ is
increasing. \label{lm1}
\end{lm}

{\bf Proof}. It suffices to prove the results after ``moreover''.
The first result that if $F\in {\cal F}(I)$ is increasing then $F^n$
is increasing can be found from \cite[Lemma 2.2]{Xu-Nk-Zhang}. The
second result is a simple corollary of the following claim: If
$F_1,F_2\in {\cal F}(I)$ are decreasing then $F_2\circ F_1$ is
increasing; if $F_1$ is decreasing and $F_2$ is increasing then
$F_1\circ F_2$ is decreasing.

Suppose that $F_1,F_2\in {\cal F}(I)$ are decreasing. For
$x_1,x_2\in I$ with $x_1<x_2$, we have ${\rm min}F_2(x_1)>{\rm
max}F_2(x_2)$ and ${\rm min}F_1({\rm max}F_2(x_2))>{\rm max}F_1({\rm
min}F_2(x_1))$ since $F_1,F_2$ are both decreasing. Clearly,
$$
{\rm min}F_1({\rm max}F_2(x_2))={\rm min}F_1\circ F_2(x_2),
$$
$$
{\rm max}F_1({\rm min}F_2(x_1))={\rm max}F_1\circ F_2(x_1),
$$
implying that ${\rm min}F_1\circ F_2(x_2)>{\rm max}F_1\circ F_2(x_1)$. Thus, $F_1\circ F_2$ is increasing. In order to prove the second result of the claim, consider the case that $F_1$ is decreasing and $F_2$ is increasing.
For $ x_1,x_2\in I$ with $x_1<x_2$, we have ${\rm max} F_1(x_2)<{\rm min}F_1(x_1)$ and
${\rm max} F_2(x_1)<{\rm min}F_2(x_2)$. It follows that
$
{\rm max}F_1({\rm min}F_2(x_2))<{\rm min} F_1({\rm max}F_2(x_1)).
$
Consequently, ${\rm max}F_1\circ F_2(x_2)<{\rm min}F_1\circ F_2(x_1)$, implying that
$F_1\circ F_2$ is decreasing. The proof is completed. $\hspace{\stretch{1}}\Box$

For each multifunction $F\in {\cal F}(I)$, let $J(F)$ denote the set of jumps of $F$.

\begin{lm}
Let $F,G\in {\cal F}(I)$ and $x\in I$.
Then $x\in J(G\circ F)$ if and only if either $x\in J(F)$ or $F(x)$ is a singleton contained in $J(G)$.
\label{lm7}
\end{lm}

{\bf Proof}. Suppose that $x_0\in J(G\circ F)$, i.e., $\#G(F(x_0))\geq 2$. By the monotonicity of $G$, either $x_0\in J(F)$ or
$F(x_0)\cap J(G)\neq \emptyset$.

Assume that $x\in J(F)$. By the monotonicity we get $x\in J(G\circ F)$.
Furthermore, for any $x_0\in I$ such that $F(x_0)\cap J(G)\neq \emptyset$,
it implies that $y_0\in J(G)$ and $y_0\in F(x_0)$ for some $y_0\in I$.
Therefore, $\#G(F(x_0))\geq 2$,
and thus $x_0\in J(G\circ F)$. This completes the first assertion and the second
one follows immediately.
$\hspace{\stretch{1}}\Box$

Lemma~\ref{lm7} implies that for any $F\in {\cal F}(I)$ the sequence
$(J(F^n)_{n\in \mathbb{N}_0})$ is increasing,
i.e.,
$J(F^n)\subset J(F^{n+1})$ for all $n\in \mathbb{N}_0$.
In what follows,
we only consider functions $F\in {\cal F}(I)$ having
at most finitely many jumps. Then $(\#J(F^n)_{n\in \mathbb{N}_0})$
is an increasing sequence of nonnegative integers.
If $\#J(F^k)=\#J(F^{k+1})$ for some integers $k\in \mathbb{N}_0$,
we call the least integer $k$ the {\it intensity} of $F$ and use $\zeta(F)$ to denote it; otherwise, we define $\zeta(F)=+\infty$.

Given a real number $x$ we denote by $[x]$ the least integer not less than $x$.

\begin{lm}
Let $F\in {\cal F}(I)$ and assume that $\zeta(F)<+\infty$. Then
$
J(F^{\zeta(F)})=J(F^{\zeta(F)+i})
$
and
$
\zeta(F^i)=[\zeta(F)/i]
$
for every $i\in \mathbb{N}$.
\label{lm2}
\end{lm}

{\bf Proof}.
To show the first property it suffices to prove that if
$J(F^{k+1})=J(F^k)$ for an integer $k\in \mathbb{N}_0$ then $J(F^{k+2})\subset J(F^{k+1})$.
In fact,
for a reduction to absurdity, let
$x\in J(F^{k+2})\backslash J(F^{k+1})\ne \emptyset$.
Note that
$x\in J(F\circ F^{k+1})$ but $x\not\in J(F^{k+1})$. By Lemma~\ref{lm7} we have $F^{k+1}(x)\in J(F)$.
Then,
$$
\#F^{k+1}(F(x))=\#F(F^{k+1}(x))\geq 2,
$$
so $F(x)\subset J(F^{k+1})=J(F^k)$. Hence, $\#F^{k+1}(x)=\#F^k(F(x))\geq 2$, that is, $x\in J(F^{k+1})$,
a contradiction to the indirect assumption.

To verify the second equality fix an $i\in \mathbb{N}$ and put $r=[\zeta(F)/i]$. Since
$$
i(r-1)<\zeta(F)\leq ir<i(r+1),
$$
we have
$
J(F^{i(r-1)})\subsetneq J(F^{\zeta(F)})=J(F^{ir})=J(F^{i(r+1)}),
$
that is,
$$
J((F^i)^{r-1})\subsetneq J((F^i)^r)=J((F^i)^{r+1}).
$$
This implies that $\zeta(F^i)=r=[\zeta(F)/i]$, the second result. $\hspace{\stretch{1}}\Box$




In the sequel, we concentrate to those multifunctions $F\in {\cal
F}(I)$ with intensity $\zeta(F)=1$, called {\it exclusive multifunctions}
and use ${\cal F}_1(I)$ to collect all of them.
For convenience, suppose that $F$ has the set of
jumps
$$
J(F):=\{c_1,\ldots, c_m\}
$$
such that $a=:c_0<c_1<\ldots<c_m<c_{m+1}:=b$. Moreover, let
$$
I(F):=\{I_i: i=0,1,\ldots, m\},
$$
where $I_i:=(c_i, c_{i+1})$.

\begin{lm}
Let $F\in {\cal F}_1(I)$. Then for every $i\in \{0,1,...,m\}$ there
exists $j\in \{0,1,...,m\}$ such that $F(I_i)\subset I_j$.
\label{Fij}
\end{lm}

{\bf Proof}. For an indirect proof, assume that there exists an $i\in \{0,1,..., m\}$
such that $F(I_{i})\not\subset I_{j}$ for all $j=0,1,...,m$.
It implies that there are two interior points $x_1,x_2\in I_i$
such that $F(x_1)\in I_k$ and $F(x_2)\in I_{k+1}$,
where $I_k$ and $I_{k+1}$ are two consecutive subintervals in $I(F)$
because of the continuity of $F$.
Let $\{c\}=cl(I_{k})\cap cl(I_{k+1})$, i.e.,
the common end-point of the closure of two subintervals.
Clearly, $c\in J(F)$.
Then by the continuity of $F$ on $I_i$,
there exists an $x_3$ between $x_1$ and $x_2$,
which is surely a continuous point of $F$ such that
\begin{eqnarray}
F(x_3)=c.
\label{Fdc}
\end{eqnarray}
By Lemma~\ref{lm7},
$J(F^2)=J(F)\cup \{x\in I\mid F(x)\cap J(F)\neq \emptyset\}$.
It follows that $x_3\in J(F^2)$ because (\ref{Fdc}) implies that
$x_3\in \{x\in I\mid F(x)\cap J(F)\neq \emptyset\}$. On the other hand,
$x_3\not\in J(F)$. It means that
$\#J(F)<\#J(F^2)$,
a contradiction to the fact that $\zeta(F)= 1$. $\hspace{\stretch{1}}\Box$

Remark that the set $J(F)$ partitions the interval $I$ into $m+1$ subintervals and
$F$ is single-valued, strictly monotone and continuous on each $I_i$.

\begin{lm}
Suppose that $F\in {\cal F}_1(I)$ is strictly increasing. Then there
is a subinterval $I'\in I(F)$ such that $F(I')\subset I'$.
\label{lmFII}
\end{lm}

{\bf Proof}.
By Lemma~\ref{Fij},
for every $i\in \{0,1,...,m\}$
there exists an $j\in \{0,1,...,m\}$ such that $F(I_i)\subset I_j$.
For an indirect proof, we assume that the corresponding integer
$j\not= i$
for every $i$.
Lemma~\ref{Fij} implies that there is a sequence $\{i_1,...,i_{k},...\}$, where $i_k$'s belong to $\{0,1,..., m\}$ such that
\begin{eqnarray*}
I_i \xrightarrow[]{F} I_{i_1} \xrightarrow[]{F} ...
\xrightarrow[]{F} I_{i_{k-1}} \xrightarrow[]{F} I_{i_{k}}
\xrightarrow[]{F} ...,
\label{IF}
\end{eqnarray*}
where $J' \xrightarrow[]{F} J''$ means $F(J')\subset J''$.
Since $\# I(F)<\infty$, the sequence contains a finite sub-sequence $\mathcal{N}:=\{i_{k_0}, i_{k_1},..., i_{k_\tau}\}$
of distinct integers, $1\le \tau\le m-1$,
 such that
\begin{eqnarray}
I_{i_{k_0}} \xrightarrow[]{F} I_{i_{k_1}} \xrightarrow[]{F} I_{i_{k_2}}
\xrightarrow[]{F}... \xrightarrow[]{F} I_{i_{k_\tau}}
\xrightarrow[]{F} I_{i_{k_0}}.
\label{IFj}
\end{eqnarray}
For each $I_\alpha\in I(F)$ let $\iota(I_\alpha)$ denote its index, i.e., $\iota(I_\alpha)=\alpha$.
Then we reset integers in sequence $\mathcal{N}$ and recall them as
$z_1,z_2,...,z_\tau,z_{\tau+1}$ satisfying $z_1<z_2<...<z_\tau<z_{\tau+1}$.
It follows from (\ref{IFj}) and the definition of index that
$$
\iota(I_{z_1})<\iota(F(I_{z_1}))~~~\mbox{ and }~~~\iota(F(I_{z_{\tau+1}}))<\iota( I_{z_{\tau+1}}).
$$
This enables
us to see that
$\iota(I_{\ell})< \iota(F(I_{\ell}))$ and $\iota(F(I_{\ell+1}))<\iota( I_{\ell+1})$
for an integer $\ell\in \{z_1,z_2,...,z_\tau,z_{\tau+1}\}$,
implying that
\begin{eqnarray*}
F(y)<c_{\ell}<c_{\ell+1}<F(x),~~~~\forall~ x\in I_{\ell},~~ y\in I_{\ell+1}.
\label{FccF}
\end{eqnarray*}
This contradicts to the assumption that $F$ is strictly increasing. Therefore the proof is completed.
$\hspace{\stretch{1}}\Box$

Let $\Lambda(F)$ denote the set of all open subintervals $J\in I(F)$
satisfying that $F(J)\subset J$, where $F\in {\cal F}_1(I)$ is
strictly increasing. By the proof of Lemma~\ref{lmFII}, for every
$I_i\in I(F)\backslash \Lambda(F)$, there is a subinterval $I_j\in
\Lambda(F)$ such that $I_i$ is mapped into $I_j$ by $F$-actions.
This fact shows a correspondence
from the set $I(F)$ of all subintervals to $\Lambda(F)$,
as indicated in the following lemma.

\begin{lm}
Each strictly increasing $F\in {\cal F}_1(I)$ determines a
correspondence $\kappa_F: \{0,1,..., m\}\to \{1,...,m\}$ such that
for every $i\in \{0,1,..., m\}$ there exists $K\in \Lambda(F)$
satisfying {\rm (i)} $F^{\kappa_F(i)}(I_i)\subset K$, {\rm
(ii)} $F^{k}(I_i)\not\subset K$ for $1\le k\le \kappa_F(i)-1$, and
{\rm (iii)} $F^{p}(I_i)\cap F^{q}(I_i)=\emptyset$ for $1\le
p\not=q\le \kappa_F(i)$. \label{lmFik}
\end{lm}

{\bf Proof}. Since $\# I(F)<\infty$, from the proof of
Lemma~\ref{lmFII} we see that there exists a finite
sub-sequence $\{i_1,...,i_k\}$ of $\{0,1,...,m\}$
such that
\begin{eqnarray}
&&I_{i} \xrightarrow[]{F} I_{i_{1}} \xrightarrow[]{F} I_{i_{2}}
\xrightarrow[]{F}... \xrightarrow[]{F} I_{i_{k}}
\xrightarrow[]{F} I_{i}:=I_{i_{k+1}},
\label{ij}
\\
{\rm or}~~~
&&
I_{i} \xrightarrow[]{F} I_{i_{1}} \xrightarrow[]{F} I_{i_{2}}
\xrightarrow[]{F}... \xrightarrow[]{F} I_{i_{k}}
\xrightarrow[]{F} I_{i_{k}},
\label{i}
\end{eqnarray}
where $1\le k\le m-1$. Obviously, cases (i)-(iii) are true for
$I_i\in \Lambda(F)$. If $I_{i_{1}}= I_i$, it follows that $I_i=K$ and we define the correspondence $\kappa_F: i\mapsto 1$. It is easy to see that $\kappa_F$ satisfies results (i)-(iii). If $I_{i_{1}}\neq I_i$, then the sequence
$\{I_{i_{\tau}}\}_{\tau=1}^{k+1}$ is monotone by the strictly monotonicity
of $F$, which contradicts to the situation (\ref{ij}). Therefore,
situation (\ref{i}) enable us to define a correspondence $\kappa_F: i\mapsto k$ such that results (i)-(iii) hold.
$\hspace{\stretch{1}}\Box$

For convenience, the open subinterval $K\in \Lambda(F)$ indicated in
Lemma~\ref{lmFik} is called an {\it absorbing interval}
of $F$.
Usually, $F$ may have more than one absorbing intervals. Furthermore, for each $i\in \{0,1,...,m\}$ the natural number $\kappa_F(i)$,
defined in Lemma~\ref{lmFik},
is called the {\it absorbing time} of $F$ on the subinterval $I_i\in I(F)$. The natural number
$\ell(F):=\max\{\kappa_F(i): i=0,1,...,m\}$
is called the {\it absorbing time} of $F$.
The absorbing time is actually the minimal number of iteration of $F$
such that the image of the whole interval is finally covered by $\Lambda(F)$.

Note that the `characteristic interval' considered in \cite{Liu-Zhang11,
ZhangYang1983AMS,ZWN-Pol} for PM functions $F: I\rightarrow I$ (a
class of strictly piecewise monotone functions) is also a kind of
absorbing intervals but quite different from the above mentioned
one. As known in \cite{ZWN-Pol}, each PM function whose iteration
does not increase the number of forts has a unique characteristic
interval, a closed sub-interval bounded by either consecutive forts
or end-point, which covers the range of $F$, i.e., the domain $I$
is mapped into the characteristic interval by one-step iteration of
$F$. In contrast, our `absorbing interval' for multifunctions is
open, non-unique and not one-step.


\section{Increasing roots of increasing multifunctions}
\setcounter{equation}{0}

In this section we discuss on strictly increasing iterative roots of strictly increasing multifunctions $F\in {\cal F}_1(I)$.
The following lemma shows that iterative roots also preserve the partition of $F$
as $F$ does in Lemma~\ref{Fij}.


\begin{lm}
Let $f$ be a usc
iterative root of $F\in {\cal F}_1$. Then for each subinterval $I_i\in I(F)$ there is an integer
$j\in \{0,1,...,m\}$ such that $f(I_i)\subset I_j$. \label{fij}
\end{lm}

{\bf Proof}.
For an indirect proof, assume that there exists an $i\in \{0,1,..., m\}$
such that $f(I_{i})\not\subset I_{j}$ for all $j=0,1,...,m$.
It follows from Lemma ~\ref{lm7} that $J(f)\subset J(F)$, which implies that $f$ is continuous
and strictly monotone on $I_i$. Similar to the proof of Lemma~\ref{Fij},
it is easy to find a continuous point $x_4$ of $F$ in $I_i$ and a jump $c'$ of $F$ such that $ f(x_4)=c'.$
Note that $J(f\circ F)=J(F\circ f)=J(f)\cup \{x\in I\mid f(x)\cap J(F)\neq \emptyset\}$, which shows that $x_4\in J(f\circ F)$.
However, $x_4\not\in J(F)$ implies that $\#J(F)<\#J(F\circ F)$, a contradiction to the fact that $\zeta(F)= 1$.
Therefore the proof is completed.
$\hspace{\stretch{1}}\Box$

%

\begin{lm}
Let $f\in {\cal F}(I)$ be an increasing iterative root of $F\in
{\cal F}_1(I)$. Then there is a subinterval $I'\in I(F)$ such that
$f(I')\subset I'$. Furthermore, each $f$ defines a correspondence
$\kappa_f(i): \{0,1,..., m\}\to \{1,...,m\}$ such that for each
$i\in \{0,1,..., m\}$ there exists $K\in \Lambda(F)$ satisfying {\rm (i)} $f^{\kappa_f(i)}(I_i)\subset K$, {\rm (ii)}
$f^{k}(I_i)\not\subset K$ for $1\le k\le \kappa_f(i)-1$, and {\rm
(iii)} $f^{p}(I_i)\cap f^{q}(I_i)=\emptyset$ for $1\le p\not=q\le
\kappa_f(i)$. \label{lm16}
\end{lm}

{\bf Proof}. Let $f\in {\cal F}(I)$ be an increasing iterative root of $F$. By Lemma~\ref{fij}, for every $i\in \{0,1,...,m\}$
there exists $j\in \{0,1,...,m\}$ such that $f(I_i)\subset I_j$. Since $f$ is strictly increasing,
it follows from Lemma~\ref{lmFII} that there is a subinterval $I'\in I(F)$ such that $f(I')\subset I'$.
This proves the first assertion.
Furthermore, let $K$ denote the subinterval $I'$.
Obviously, $K\in \Lambda(F)$.
Since $\#I(F)<\infty$, we infer from the proof of Lemma~\ref{lmFik}
that there is a finite sub-sequence of $I(F)$
such that
(\ref{ij})-(\ref{i}) hold for $f$, which leads to the results of (i)-(iii). The proof is completed.
$\hspace{\stretch{1}}\Box$

Although the absorbing interval of a given exclusive multifunctions is not unique, we can partition the whole interval $I$ by
fixed points of $F$ and discuss iterative roots of $F$ restricted to
those partitioned subintervals separately, each of which contains a
unique absorbing interval.

\begin{lm}
Let $F\in {\cal F}_1(I)$ and $f\in {\cal F}(I)$ be an increasing
iterative root of $F$. If there exists $c\in (a,b)$ such that $c\in
F(c)$, then $c\in f(c)$. \label{lm17}
\end{lm}

{\bf Proof}. Let $f\in {\cal F}(I)$ be an increasing iterative root of $F$ of order $n$.
If $c$ is a single-valued point of $F$, i.e. $F(c)=c$, then we have $f(c)=c$
by the monotonicity of $f$.

In the case that $c$ is a set-valued point of $F$,
for an indirect proof we assume that $c\not\in f(c)$. Let $S_1:=\{x\in f(c)\mid x<c\}$ and $S_2:=\{x\in f(c)\mid x>c\}$ respectively.
Obviously, $f(c)=S_1\cup S_2$. Since $f$ is strictly increasing, we have $\max f(S_1)<c$ and  $\min f(S_2)>c$,
implying that $c\not\in f(f(c))$. We can inductively prove that
$c\not\in f^i(c)$ for all $i=1,2,...$, but it contradicts to the fact that $c\in F(c)= f^n(c)$. The proof is completed.
$\hspace{\stretch{1}}\Box$


Lemma~\ref{lm17} shows that iterative roots of $F$ satisfying $c\in F(c)$ can be found on subinterval $[a,c]$ and $[c,b]$
separately. So we suppose,
without loss of generality, that either ${\rm max}F(x)<x$ or ${\rm
min} F(x)>x$ for all $x\in(a,b)$. In what follows, it suffices to
consider $F$ under the assumption that
\begin{description}
\item[($\cal{H}$)] $F$ is strictly increasing, ${\rm max}F(x)<x$ for all $x\in (a,b)$ and $\max F(b)\le b$.
\end{description}
The other case can be converted to the same by the function $G(x):=(b+a)-F(b+a-x)$. Clearly, under ($\cal{H}$) we have $F(a)=a$.

\begin{lm}
Suppose that $F\in {\cal F}_1(I)$ and hypothesis ($\cal{H}$) holds.
Then the interval $K:=(a,c_1)\in I(F)$ is the unique absorbing
interval of $F$. Furthermore, $\Lambda(f)=\Lambda(F)=K$ if $f\in
{\cal F}(I)$ is an increasing iterative root of $F$. \label{lmfII}
\end{lm}

{\bf Proof}.
The existence of absorbing interval, denoted by $K$, is given in Lemma~\ref{lmFik}.
By Lemma~\ref{fij} and assumption (${\cal H}$), $F(K)\subset K$ if and only if $K=(a,c_1)$.
It implies the uniqueness of absorbing interval. In order to prove the second part,
noting that $f\in {\cal F}(I)$ is strictly increasing, we claim that ${\rm max}f(x)<x$ for all $x\in (a,b)$.
Otherwise, the assumption that ${\rm max}f(x_0)\geq x_0$ for some $x_0\in (a,b)$ leads to the fact
that ${\rm max}F(x_0)={\rm max}f^n(x_0)\geq x_0$, which is a contradiction.
Therefore, $\Lambda(f)=\Lambda(F)=K$ and the proof is completed. $\hspace{\stretch{1}}\Box$

By the above lemmas, we can find iterative roots of $F$
by constructing iterative roots on the absorbing interval first and then
extending them to the whole domain.

\begin{thm}\label{thm4}
Suppose that $F\in {\cal F}_1(I)$ and hypothesis ($\cal{H}$) holds.
Then every iterative root $f\in {\cal F}(I)$ of order $n$ of $F$ on
$I\setminus J(F)$ is an extension from an increasing iterative root
of $F$ of the same order on the absorbing interval $K$.
\end{thm}


{\bf Proof}. Consider the restrictions of increasing function $f\in {\cal F}(I)$
to the interval $I_i\in I(F)$ and let $f_i:=f|_{I_i}$, where $I_i\in I(F)$, and $f|_{K}:=\phi$ for $\phi$ being an increasing iterative root of $F$ of order $n$ on $K$.
By Lemma~\ref{lm16}, for each $I_i$ there
is a finite sub-sequence $\{i_1,\ldots,i_{k-1}\}$ of $\{0,1,...,m\}$ such that
$$
I_{i} \xrightarrow[]{f} I_{i_{1}} \xrightarrow[]{f} I_{i_{2}}
\xrightarrow[]{f}... \xrightarrow[]{f} I_{i_{k-1}}
\xrightarrow[]{f} K.
$$
Consequently, for every $x\in I_{i_{k-1}}$ we have
$F(x)=\phi^{n-1}\circ f_{i_{k-1}}(x)$
because $\phi(K)\subset K$. It follows that
\begin{eqnarray}\label{fik1}
f_{i_{k-1}}(x)=F|_{K}^{-1}\circ \phi\circ F|_{I_{i_{k-1}}}(x).
\label{Ff2}
\end{eqnarray}
Furthermore, for every $x\in I_{i_{k-2}}$ we have
\begin{eqnarray}
\begin{array}{lll}
f_{i_{k-2}}(x)&=&f_{i_{k-1}}^{-1}\circ \phi^{-n+2}\circ F|_{I_{i_{k-2}}}(x)
\\
\vspace{-0.2cm}
\\
&=&F|_{I_{i_{k-1}}}^{-1}\circ \phi^{-1}\circ F|_{K}\circ \phi^{-n+2}\circ F|_{I_{i_{k-2}}}(x)
\\
\vspace{-0.2cm}
\\
&=&F|_{I_{i_{k-1}}}^{-1}\circ \phi\circ F|_{I_{i_{k-2}}}(x).
\end{array}
\label{fik2}
\end{eqnarray}
One can check by induction that
\begin{eqnarray}
\left\{
\begin{array}{lll}
f_{i_j}&=&F|_{I_{i_{j+1}}}^{-1}\circ \phi\circ F|_{I_{i_{j}}}, ~~j=1,2,\ldots, k-2,
\\
\vspace{-0.2cm}
\\
f_{i}&=&F|_{I_{i_{1}}}^{-1}\circ \phi\circ F|_{I_{i}}.
\end{array}
\right.
\label{fik3}
\end{eqnarray}
Therefore, in view of (\ref{fik1})-(\ref{fik3}), the
iterative root $f_{*}:=f\mid_{I\setminus J(F)}$ is an extension of $\phi$ from $K$ in the form of
\begin{equation}\label{deff}
f_*(x)=\left\{
\begin{array}{ll}
\phi(x), \ \ & \ \ {\rm as}\ \  x\in K,
\\
F|_{I_{i_{1}}}^{-1}\circ \phi\circ F|_{I_{i}},\ \ & \ \ {\rm as}\ \  x\in I_i\in I(F)\backslash \{K\}.
\end{array}
\right.
\end{equation}
The proof is completed.
$\hspace{\stretch{1}}\Box$

The proof of Theorem~\ref{thm4} gives the construction of increasing iterative roots of $F$ on $I\setminus J(F)$,
the place without a jump. In what follows, we will consider iterative roots defined at jumps.
For each jump $c_i\in J(F)$, $i=1,2,\ldots, m$,
under hypothesis ($\cal{H}$) we need to consider the four cases
as mentioned in the Introduction, i.e.,
\begin{description}
\item{\bf (J1)}
$J(F)\cap F(c_i)=\emptyset$,

\item{\bf (J2)}
$J(F)\cap F(c_i)=\{c_i\}$,

\item{\bf (J3)}
$J(F)\cap F(c_i)=\{c_{k_1}, c_{k_2}, \ldots, c_{k_\ell}\}$,
where $c_{k_1},c_{k_2},\ldots,c_{k_\ell}\in J(F)\backslash \{c_i\}$
and $1\le \ell<m-1$ is an integer,
and

\item{\bf (J4)}
$J(F)\cap F(c_i)=
\{ c_{k_1}, c_{k_2}, \ldots, c_{k_\ell}, c_i\}$,
\end{description}

The following Theorems~\ref{thm7}, \ref{thm16}  and \ref{thm8} are devoted to cases {\bf (J1)}, {\bf (J2)} and {\bf (J3)} respectively
for iterative roots on $(I\setminus J(F))\cup \{c_i\}$.

\begin{thm}
Suppose that $F\in {\cal F}_1(I)$ with hypothesis ($\cal{H}$) lies in
case {\bf (J1)} for a given $c_i\in J(F)$ and
that
$f_*$ is defined in (\ref{deff}) and
maps $I\setminus J(F)$ onto itself. Then the
multifunction
\begin{eqnarray}\label{fdefineJ1}
f(x):=\left\{
\begin{array}{ll}
f_{*}(x), \ \ & \ \ {\rm as}\ \ x\in I\backslash J(F),
\\
f_{*}^{-(n-1)}\circ F(c_i),\ \ & \ \ {\rm as}\ \  x=c_i,
\end{array}
\right.
\end{eqnarray}
is a strictly increasing usc
iterative root
of $F$ of order $n$ on $(I\setminus J(F))\cup \{c_i\}$.
\label{thm7}
\end{thm}

\begin{thm}
Suppose that $F\in {\cal F}_1(I)$ with hypothesis ($\cal{H}$) lies in
case {\bf (J2)} for a given $c_i\in J(F)$, the set $F(c_i)$
is a compact interval,
and
$f_*$ is defined in (\ref{deff}).
Then the multifunction
\begin{eqnarray}\label{fdeinJ2}
f(x):=\left\{
\begin{array}{ll}
f_*(x), \ \ & \ \ {\rm as}\ \ x\in I\backslash J(F),
\\
\,[\lim_{s\to c_{i}-}f_*(s),c_i],\ \ & \ \ {\rm as}\ \  x=c_i,
\end{array}
\right.
\end{eqnarray}
is a strictly increasing usc
iterative root of $F$
of order $n$ on $(I\setminus J(F))\cup \{c_i\}$.
\label{thm16}
\end{thm}

\begin{thm}
Suppose that $F\in {\cal F}_1(I)$ with hypothesis ($\cal{H}$) lies in
case {\bf (J3)} for a given $c_i\in J(F)$. Then $F$ has no strictly
increasing usc
iterative roots of order
$n>m-\ell+1$ on $I$. \label{thm8}
\end{thm}

{\bf Proof of Theorem~\ref{thm7}}.
According to formula (\ref{fdefineJ1}), the strictly monotonicity of $f$ on $I\backslash J(F)$ is obvious by Theorem~\ref{thm4}. Moreover,
$$
f^n(c_i)=f_*^{n-1}\circ f_*^{-(n-1)}\circ F(c_i)=F(c_i)
$$
since $f_*$ maps $I\setminus J(F)$ onto itself. Hence, the function $f$ defined in (\ref{fdefineJ1}) is an iterative root of $F$ of order $n$ on $(I\setminus J(F))\cup \{c_i\}$.
As indicated in \cite[Lemma 2]{Li-Jar-Zhang09Publ}, a strictly increasing multifunction $F$ is
usc
at $c_i$
if and only if
$F(c_i-)$ (resp. $F(c_i+)$) is the smallest (resp. greatest) element of $F(c_i)$.
Then, by the monotonicity of $\phi$ and $F$, we obtain
\begin{eqnarray*}
&&\inf f(c_i)=\inf f_*^{-(n-1)}\circ F(c_i)=\lim_{x\to c_{i}-} f_*^{-(n-1)}\circ F(x)=\lim_{x\to c_{i}-} f(x),
\\
&&\sup f(c_i)=\sup f_*^{-(n-1)}\circ F(c_i)=\lim_{x\to c_{i}+} f_*^{-(n-1)}\circ F(x)=\lim_{x\to c_{i}+} f(x).
\end{eqnarray*}
Therefore, $f$ is increasing and usc
at $c_i$. The proof is completed.
$\hspace{\stretch{1}}\Box$


{\bf Proof of Theorem~\ref{thm16}}. Suppose that $F\in {\cal F}_1(I)$
with hypothesis ($\cal{H}$) lies in case {\bf (J2)} for a given
$c_i\in J(F)$. Then $c_i$ is the right endpoint of $I$. Since the
proof of Theorem~\ref{thm4} gives the construction of increasing
iterative roots of $F$ on $I\setminus J(F)$, one can define $f=f_*$
on $I\backslash J(F)$, where $f_*$ is given in (\ref{deff}).

Note that the set $F(c_i)$ is a compact interval, i.e., $F(c_i)=[\inf F(c_i), \sup F(c_i)]$. Then $F(c_i)=[\lim_{x\to c_{i}-}F(x), c_i]$
because $F$ is usc
at $c_i$ satisfying {\bf (J2)}. Let
$$
\Omega:=[\lim_{x\to c_{i}-}f_*(x),c_i].
$$
By hypothesis ($\cal{H}$) and ${\bf (J2)}$ we get
\begin{eqnarray*}
\label{Mi1Fh}
c_{i-1}<\lim_{x\to c_{i}-}F(x)<\lim_{x\to c_{i}-}f_*(x)<c_i
\end{eqnarray*}
because $f_*$ is strictly increasing, it follows that $\Omega\ne \emptyset$. Therefore, by \cite[Lemma 2]{Li-Jar-Zhang09Publ}, formula (\ref{fdeinJ2})
defines
a strictly increasing multifunction $f$ on $(I\setminus J(F))\cup \{c_i\}$,
which is usc.
Further,
one can check that
$f^n(x)=f_*^n(x)=F(x)$ for all $x\in I\setminus J(F)$ by Theorem~\ref{thm4} and
\begin{eqnarray}
\begin{array}{lll}
f^n(c_i)&=&f^{n-1}(\Omega)
=\bigcup_{j=1}^{n-1}(f_*^{j}([\lim_{x\to c_{i}-}f_*(x),c_i)))\cup\Omega
\\
\vspace{-0.2cm}
\\
&=&[\lim_{x\to c_{i}-}f^n_*(x),\lim_{x\to c_{i}-}f^{n-1}_*(x))\cup\ldots\cup[\lim_{x\to c_{i}-}f_*(x),c_i]
\\
\vspace{-0.2cm}
\\
&=&[\lim_{x\to c_{i}-}f^n_*(x),c_i]
=[\lim_{x\to c_{i}-}F(x),c_i]=F(c_i),
\end{array}
\label{revisedomega}
\end{eqnarray}
which follows that $f$ defined in (\ref{fdeinJ2}) is an $n$-th order iterative root of $F$ on $(I\setminus J(F))\cup \{c_i\}$.
This completes the proof.
$\hspace{\stretch{1}}\Box$


{\bf Proof of Theorem~\ref{thm8}}.
Assume that $f\in {\cal F}(I)$ is an increasing iterative root of $F$ of order $n$. Let
\begin{eqnarray}\label{sf}
S_j:=f^j(c_i)\cap J(F), \quad ~ j=1,2,\ldots,n,
\end{eqnarray}
none of which is empty by Lemma~\ref{fij} because $f^n(c_i)=F(c_i)\cap J(F)\neq \emptyset$.
Under the assumption of case {\bf (J3)}, we claim that
\begin{eqnarray}
d_j\not\in f(d_j),\quad \forall d_j\in S_j.
\label{3.2}
\end{eqnarray}
In fact,
(\ref{3.2}) is true obviously for $j=n$; otherwise, we have
$d_n\in f(d_n)$, implying that $d_n\in F(d_n)$,
a contradiction to the fact that $d_n\in \{c_{k_1}, c_{k_2}, \ldots, c_{k_l}\}\subset F(c_i)$
because $F$ is strictly monotone.
Further, for an indirect proof, we assume that
there exists an integer $j^*\in \{1,2,\ldots,n-1\}$ such that
\begin{eqnarray}
d_{j^*}\in f(d_{j^*}).
\label{dj}
\end{eqnarray}
It follows by iteration that $d_{j^*}\in F(d_{j^*})$.
On the other hand,
$d_{j^*}\in S_{j^*}=f^{j^*}(c_i)\cap J(F)$ implies that $d_{j^*}\in f^n(c_i)=F(c_i)$ because we have (\ref{dj}) and $j^*<n$.
Thus,
$d_{j^*}\in F(d_{j^*})\cap F(c_i)$,
a contradiction to the monotonicity of $F$ since $d_{j^*}\neq c_i$.
Therefore, the claimed (\ref{3.2}) is proved.
Since $f$ is strictly increasing, we infer from (\ref{3.2}) that $S_i\backslash S_j\neq\emptyset$ for any $i\neq j\in\{1,2,\ldots,n\}$.
Then,
by the fact that $\#J(F)=m$ and $\#S_n=\ell$, we get $m\geq \ell+n-1$. This completes the proof.
$\hspace{\stretch{1}}\Box$


Up to now, the problem of iterative roots is not solved yet in  case {\bf (J4)}.


\section{Decreasing roots of increasing multifunctions}

In this section, we discuss the decreasing iterative roots of
strictly increasing multifunctions $F\in {\cal
F}_1(I)$. As known in \cite[Theorem
11.2.5]{KCG1990}, a strictly increasing continuous self-mapping has
no decreasing iterative roots of odd order but may have roots of
even order. This fact is also true for strictly monotone
semi-continuous multifunctions by Lemma~\ref{lm1}. Let $f\in {\cal
F}(I)$ be a decreasing iterative root of $F\in {\cal F}_1(I)$ of order
$n=2k\ge 2$. Then $g:=f^2$ is a strictly increasing and upper
semi-continuous iterative root of $F$, which was discussed in
section 3. Therefore, we shall confine ourselves to consider the functional equation
\begin{eqnarray}
f^2(x)=F(x),~~~\forall x\in I,
\label{square root}
\end{eqnarray}
where $F\in {\cal F}_1(I)$ is a given strictly increasing
multifunction and decreasing multifunction $f\in {\cal F}(I)$ is
unknown.

There are more difficulties in the case of decreasing roots because
$f$ may have no absorbing intervals anymore. For example, consider
$$
F(x):=\left\{
\begin{array}{ll}
\frac{1}{4}x+\frac{1}{8}, \ \ & \ \ {\rm as}\ \  x\in [0,\frac{1}{2}),
\\
\,[\frac{1}{4},\frac{3}{4}],\ \ & \ \ {\rm as}\ \  x = \frac{1}{2},
\\
\frac{1}{4}x+\frac{5}{8}, \ \ & \ \ {\rm as}\ \  x\in (\frac{1}{2},1],
\end{array}
\right.
~~~
f(x):=\left\{
\begin{array}{ll}
-x+1, \ \ & \ \ {\rm as}\ \  x\in
[0,\frac{1}{2}),
\\
\,[\frac{1}{4},\frac{1}{2}],\ \ & \ \ {\rm as}\ \  x = \frac{1}{2},
\\
-\frac{1}{4}x+\frac{3}{8}, \ \ & \ \ {\rm as}\ \  x\in
(\frac{1}{2},1].
\end{array}
\right.
$$
One can check that $f\in {\cal F}(I)$ and $f: I\rightarrow 2^{I}$
is a square decreasing iterative root of $F$.
Note
$f(I_i)\not\subset I_i$ for each $I_i\in I(F)$, implying that $f$ has no absorbing intervals.

In contrast to Lemma~\ref{lm16},  we have the
following properties for decreasing roots.

\begin{lm}
Let $f\in {\cal F}(I)$ be a decreasing solution of equation
(\ref{square root}) with increasing $F\in {\cal F}_1(I)$. Then
$f(\Lambda(F))\subset \Lambda(F)$. More concretely, for each
subinterval $I'\in \Lambda(F)$ there exists $J'\in \Lambda(F)$ such
that $f(I')\subset J'$ and $f(J')\subset I'$. In particular, $I'=J'$
if $f$ is a self-mapping on $I'$. \label{d5.1}
\end{lm}

{\bf Proof}.
From the definition of $\Lambda(F)$ we see that $F(I')\subset I'$ for every  $I'\in \Lambda(F)$.
Then by Lemma~\ref{fij},
$f(I')\subset J'$ for a certain subinterval $J'\in I(F)$,
which implies that $F(I')\subset f(J')$ and $f(J')\subset I'$. Thus,
\begin{eqnarray*}
F(J')=f(f(J'))\subset f(I')\subset J'
\end{eqnarray*}
by (\ref{square root}). This gives the fact that $J'\in \Lambda(F)$. Since $I'\in \Lambda(F)$ is chosen arbitrarily,
we have $f(\Lambda(F))\subset \Lambda(F)$. The proof is completed. $\hspace{\stretch{1}}\Box$

\begin{lm}
Let $f\in {\cal F}(I)$ be a decreasing solution of equation
(\ref{square root}) with increasing $F\in {\cal F}_1(I)$. Then there
exists a correspondence $\kappa^{'}_f: \{0,1,..., m\}\to
\{1,...,m\}$ such that for every $i\in \{0,1,..., m\}$ we have {\rm
(i)} $f^{\kappa^{'}_f(i)}(I_i)\in \Lambda(F)$, {\rm (ii)}
$f^{k}(I_i)\not\in \Lambda(F)$ for all $1\le k\le
\kappa^{'}_f(i)-1$, and {\rm (iii)} $f^{p}(I_i)\cap
f^{q}(I_i)=\emptyset$ for $1\le p\not=q\le \kappa^{'}_f(i)$.
\label{lm14}
\end{lm}

Similar to $\kappa_F(i)$ defined in Lemma~\ref{lmFik}, the notation $\kappa^{'}_f(i)$ defined above
gives the least number of $f$-actions for subinterval $I_i$ to be mapped by $f$ iteratively into $\Lambda(F)$.

By Lemma~\ref{lm14}, square iterative roots of $F$ can be obtained in the procedure:
{\bf Step 1:}
find  strictly decreasing iterative roots on the subintervals in $\Lambda(F)$, where $F$ is a strictly increasing self-mapping;
{\bf Step 2:} extend those roots to the whole domain.
Step 1 can be completed by \cite[Theorem 15.10]{Kuczma68} or \cite[Remark 11.2.3]{KCG1990},
where a necessary condition that $F|_{\Lambda(F)}$ is a reversing correspondence
needs to be considered.
The following Theorems~\ref{thm3}-\ref{thm10} are devoted to Step 2.

\begin{thm}
Suppose that $F\in {\cal F}_1(I)$ is strictly increasing. Then every
square decreasing iterative root $f\in {\cal F}(I)$ of $F$ on
$I\backslash J(F)$ is an extension from a strictly decreasing square
iterative root of $F$ on subintervals in $\Lambda(F)$. \label{thm3}
\end{thm}

{\bf Proof}.
For  a given square decreasing iterative root $f\in {\cal F}(I)$ of $F$,
we use the notations
$\psi:=f|_{\Lambda(F)}$ and
$f_i:=f|_{I_i}$, where $I_i\in I(F)$.
By Lemma~\ref{lm14}, for each $I_i\in I(F)$ there is a finite sequence $\{i_1,\ldots,i_{k-1}\}$ in $\{0,1,...,m\}$ such that
$$
I_{i} \xrightarrow[]{f} I_{i_{1}} \xrightarrow[]{f} I_{i_{2}}
\xrightarrow[]{f}... \xrightarrow[]{f} I_{i_{k-1}}
\xrightarrow[]{f} I'\in \Lambda(F).
$$
Consequently, for every $x\in I_{i_{k-1}}$ we have
$F(x)=\psi\circ f_{i_{k-1}}(x)$.
It follows
that
\begin{eqnarray*}
f_{i_{k-1}}(x)=\psi^{-1}\circ F|_{I_{i_{k-1}}}(x)=F|_{I'}^{-1}\circ \psi\circ F|_{I_{i_{k-1}}}(x)
\end{eqnarray*}
since $\psi$ is a square iterative root of $F$ on $I'$. One can check by induction
that
\begin{eqnarray*}
\left\{
\begin{array}{lll}
f_{i_j}&=&F|_{I_{i_{j+1}}}^{-1}\circ\psi\circ F|_{I_{i_{j}}}, ~~j=1,2,\ldots, k-2,
\\
\vspace{-0.2cm}
\\
f_{i}&=&F|_{I_{i_{1}}}^{-1}\circ\psi\circ F|_{I_{i}}.
\end{array}
\right.
\end{eqnarray*}
Therefore, the
iterative root $f^{*}:=f\mid_{I\setminus J(F)}$ is an extension of $\psi$ from the subintervals in $\Lambda(F)$ in the form of
\begin{equation}\label{deffd}
f^*(x)=\left\{
\begin{array}{ll}
\psi(x), \ \ & \ \ {\rm as}\ \  x\in \Lambda(F),
\\
F|_{I_{i_{1}}}^{-1}\circ\psi\circ F|_{I_{i}},\ \ & \ \ {\rm as}\ \  x\in I_i\in I(F)\backslash \Lambda(F).
\end{array}
\right.
\end{equation}
This completes the proof.
$\hspace{\stretch{1}}\Box$

The proof of Theorem~\ref{thm3} gives the construction of decreasing square iterative roots of $F$
on $I\setminus J(F)$, the place without a jump.
In order to define those roots at a given jump $c_i\in J(F)$, $i=1,2,\ldots, m$,
we need to consider the same four cases {\bf (J1)}-{\bf (J4)} as listed in the Introduction and section 3, respectively.

\begin{thm}
Suppose that $F\in {\cal F}_1(I)$ is strictly increasing and lies in case {\bf (J1)} for a given $c_i\in J(F)$
and that
$f^*$ is defined in (\ref{deffd}) and maps $I\setminus J(F)$ onto itself.
Then the
multifunction
\begin{eqnarray}\label{fdefineD}
f(x):=\left\{
\begin{array}{ll}
f^{*}(x), \ \ & \ \ {\rm as}\ \ x\in I\backslash J(F),
\\
(f^{*})^{-1}\circ F(c_i),\ \ & \ \ {\rm as}\ \  x=c_i,
\end{array}
\right.
\end{eqnarray}
is a strictly decreasing usc
square iterative root
of $F$ on $(I\setminus J(F))\cup \{c_i\}$.
\label{thmD1}
\end{thm}

\begin{thm}
Suppose that $F\in {\cal F}_1(I)$ is strictly increasing and lies in
case {\bf (J2)} for a given $c_i\in J(F)$,
the set $F(c_i)$ is a compact interval,
and
$f^*$ is defined in (\ref{deffd}).
Then the multifunction
\begin{eqnarray}\label{fdeind}
f(x):=\left\{
\begin{array}{ll}
f^*(x), \ \ & \ \ {\rm as}\ \ x\in I\backslash J(F),
\\
\,[\lim_{s\to c_{i}+} f^*(s),\lim_{s\to c_{i}-} f^*(s)],\ \ & \ \ {\rm as}\ \  x=c_i,
\end{array}
\right.
\end{eqnarray}
is a strictly decreasing usc
square iterative root of $F$ on $(I\setminus J(F))\cup \{c_i\}$
if
\begin{eqnarray}\label{Jf^*}
J(F)\cap (\lim_{x\to c_{i}+} f^*(x),\lim_{x\to c_{i}-} f^*(x))=c_i.
\end{eqnarray}
\label{thm10}
\end{thm}

Theorem~\ref{thm10} shows that for decreasing solutions we need
condition (\ref{Jf^*}), which means that the set-value of $f$ at
$c_i$ covers exactly one jump point, which is $c_i$ itself. We also
call such a $c_i$ a fixed point of $f$ in the inclusion sense.
Moreover, the fact that the set-value of $f$ at $c_i$ contains none
of other jumps of $F$ guarantees the independence of $c_i$ under
iteration.

{\bf Proof of Theorem~\ref{thmD1}}. Since the construction of
decreasing square iterative roots of $F$ on $I\setminus J(F)$ is
given in the proof of Theorem~\ref{thm3}, one can defined $f=f^*$ on
$I\setminus J(F)$ as (\ref{fdefineD}). The definition of
(\ref{fdefineD}) at $x=c_i$ is similar to the corresponding one
given in Theorem~\ref{thm7} with $n=2$ because the fact $J(F)\cap
F(c_i)=\emptyset$ implies that $F(c_i)\subset I\setminus J(F)$ and
therefore is included in the range of $f^*$.
$\hspace{\stretch{1}}\Box$

{\bf Proof of Theorem~\ref{thm10}}. One can define $f=f^*$ on
$I\setminus J(F)$ as (\ref{fdeind}) since the construction of
strictly decreasing square iterative roots of $F$ on $I\setminus
J(F)$ is given in the proof of Theorem~\ref{thm3}.

In addition, at $c_i$ we note that the set $F(c_i)$ is a compact interval. Then $F(c_i)=[\lim_{x\to c_{i}-}F(x), c_i]$
because $F$ is usc
at $c_i$ satisfying {\bf (J2)}. Let
$$
\Omega:=[\lim_{x\to c_{i}+}f^*(x),\lim_{x\to c_{i}-}f^*(x)].
$$
Then
$\Omega\ne \emptyset$
since $f^*$ is strictly decreasing. By \cite[Lemma 2]{Li-Jar-Zhang09Publ}, formula (\ref{fdeind})
defines
a strictly decreasing multifunction $f$ on $(I\setminus J(F))\cup \{c_i\}$,
which is usc.
Furthermore, by Theorem~\ref{thm3} we have $f^2(x)=(f^{*})^2(x)=F(x)$
 if $x\in I\setminus J(F)$.
On the other hand, for $x=c_i$,
\begin{eqnarray}
f^2(c_i)&=&f(\Omega)
=f([\lim_{x\to c_{i}+}f^*(x),c_i))\cup f(c_i)\cup f((c_i,\lim_{x\to c_{i}-}f^*(x)])
\nonumber\\
&=&[\lim_{x\to c_{i}-}f^*(x),\lim_{x\to c_{i}+}(f^*)^2(x))\cup[\lim_{x\to c_{i}+}f^*(x),\lim_{x\to c_{i}-}f^*(x)]
\nonumber\\
& &\cup [\lim_{x\to c_{i}-}(f^*)^2(x),\lim_{x\to c_{i}+}f^*(x))
\nonumber\\
&=&[\lim_{x\to c_{i}-}(f^*)^2(x),\lim_{x\to c_{i}+}(f^*)^2(x)]
\nonumber\\
&=&[\lim_{x\to c_{i}-}F(x),\lim_{x\to c_{i}+}F(x)]=F(c_i)
\label{derevisedomega}
\end{eqnarray}
by the assumption (\ref{Jf^*}).
It follows that $f$ defined in (\ref{fdeind}) is a square iterative root of $F$ on $(I\setminus J(F))\cup \{c_i\}$.
This completes the proof.  $\hspace{\stretch{1}}\Box$


\section{Iterative roots of decreasing multifunctions}


We previously discussed iterative roots (more concretely, increasing
iterative roots in section 3
and decreasing iterative roots in
section 4) for strictly increasing exclusive multifunctions, and will take a concern to iterative roots for strictly
decreasing case. It is known in the single-valued case that a
strictly decreasing continuous self-mapping has no continuous roots
of even order (see \cite[p.425-426]{KCG1990}) but may have
decreasing roots of odd order (see Theorem 11.2.4 in
\cite{KCG1990}). Is there any similar result in our considered
set-valued cases?

By Lemma~\ref{lm1}, the first single-valued result with even order
is also true for strictly decreasing usc
multifunctions, i.e., each strictly decreasing usc
multifunction has no usc
iterative roots of even order.
In what follows, we concentrate on the iterative roots for
those strictly decreasing multifunctions, where their intensities are
equal to $1$. We first prove a result on nonexistence of square
roots which are continuous on $I\backslash J(F)$.

\begin{thm}\label{thm11}
Each strictly decreasing multifunction $F\in {\cal F}_1(I)$ has no
square iterative roots which are continuous on $I\backslash J(F)$.
\end{thm}

{\bf Proof}.
For an indirect proof,  suppose that $F$ has a square iterative root
$f:I\to 2^{I}$ which is continuous on $I\backslash J(F)$.
We first claim that $f$ is single-valued on $I\backslash J(F)$.
Otherwise, there is a point $x_0\in I\backslash J(F)$ such that $\#f(x_0)\geq 2$. Given two points $p,q\in f(x_0)$, since
$F(x_0)$ is a singleton, we have $F(x_0)=f(p)=f(q)$ implying that
$F(p)=F(q)$, which contradicts to the fact that $F$ is strictly decreasing. By the monotonicity of $F$, there exists a unique point $c\in I$ such that (i) $F(c)=c$ or (ii) $c\in F(c)$. Let $F_1, F_2$
denote the restrictions of $F$ to the subintervals $[a,c)$ and
$(c,b]$ respectively. Obviously,
\begin{equation}\label{F}
F_1(x)>c ~~~~{\rm and}~~~~ F_2(x)<c
\end{equation}

For case (i), $c$ is also a fixed point of $f$ since $f(c)$ is a
singleton. Moreover, there
exists $\delta>0$ such that $f$ is continuous on $[c-\delta, c+\delta]\subset I$, and
\begin{eqnarray}
F(x)> c ~~{\rm for}~x\in [c-\delta,c), \label{>c}
\\
F(x)<c~~ {\rm for} ~x\in (c,c+\delta]. \label{<c}
\end{eqnarray}
By the continuity of $f$ on $I\backslash J(F)$, choosing
sufficient small $\delta_1$, $0<\delta_1<\delta$ such that $|f(x)-c|<
\delta$ for every $x\in [c-\delta_1, c+\delta_1]$. Let
$M_f=$max$\{f(x): x\in [c-\delta_1,c]\}$. If $M_f\leq c$, then
choose small enough $\delta_2$, $0<\delta_2<\delta_1<\delta$ such that $F(x)=f^2(x)\leq M_f\leq
c$ for every $x\in [c-\delta_2,c]$, which contradicts to
(\ref{>c}). If $M_f> c$, it implies that $f(x)>c$ for every
$x\in (c, M_f]$. Otherwise, there is a point $x_1\in (c,M_f]$ such that $f(x_1)\leq c$. It follows that
there exists a point $x_2\in [c-\delta_1,c)$ satisfying $f(x_2)=x_1$. Hence, $F(x_2)=f(x_1)\leq c$, which contradicts to (\ref{>c}).
Therefore, we can choose sufficient small $\delta_3$ that
$F(x)=f^2(x)>c$ for every $x\in(c, c+\delta_3]$, however, it
is a contradiction to (\ref{<c}).

For case (ii), i.e., $c\in F(c)$, it implies that $c\in J(F)$, say $c:=c_i$.
We first claim for every $x\in I\backslash \{c_i\}$ that $c_i\not\in f(x)$. Otherwise, assume there exists an point
$y_0\in I\backslash \{c_i\}$ such that $c_i\in f(y_0)$, then
$c_i\in F(c_i)\subset F(f(y_0))=f(F(y_0))$, which gives $y_0\in F(y_0)$ by the monotonicity of $f$.
But it is a contradiction to the uniqueness of $c_i$ that $c_i\in F(c_i)$. Therefore, by the upper semi-continuity of $F$ on $I$, $F_1(x)>c_i, F_2(x)>c_i$, or
$F_1(x)>c_i, F_2(x)<c_i$, or $F_1(x)<c_i, F_2(x)>c_i$, or $F_1(x)<c_i, F_2(x)<c_i$.
However, from (\ref{F}) we see that none of these cases holds. The proof is completed.
$\hspace{\stretch{1}}\Box$

Continuing the above answer to iterative roots of even order, we
consider the second question: Does a decreasing exclusive multifunctions have a decreasing iterative root of odd order? Since
every decreasing multifunction $F\in {\cal F}_1(I)$ is a square
iterative root of the increasing multifunction $F^2\in {\cal F}_1(I)$,
every iterative root of $F$ of odd order is also a root of $F^2$ of
even order. Note that $I(F^2)=I(F)=\{I_i: i=0,1,...,m\}$ because
$\zeta(F)=1$. Thus, similarly to Lemmas~\ref{d5.1}-\ref{lm14}, we
have the following results.

\begin{lm}
\label{zxm} Suppose that the strictly decreasing multifunction $F\in
{\cal F}_1(I)$ has a decreasing iterative root $f\in {\cal F}(I)$ of
odd order $k$. Then $f(\Lambda(F^2))\subset \Lambda(F^2)$.
\end{lm}

{\bf Proof}.
From the definition of $\Lambda(F^2)$ we see that $F^2(I')\subset I'$ for every  $I'\in \Lambda(F^2)$.
Then, by Lemma~\ref{fij},
$f(I')\subset J'$ for a certain subinterval $J'\in I(F^2)$,
which implies that $F^2(I')=f^{2k-1}\circ f(I')\subset f^{2k-1}(J')\subset I'$. Thus,
\begin{eqnarray*}
F^2(J')=f\circ f^{2k-1}(J')\subset f(I')\subset J'.
\end{eqnarray*}
This gives the fact that $J'\in \Lambda(F^2)$. Since $I'\in \Lambda(F^2)$ is chosen arbitrarily,
we have $f(\Lambda(F^2))\subset \Lambda(F^2)$. The proof is completed. $\hspace{\stretch{1}}\Box$

\begin{lm}
Suppose that the strictly decreasing multifunction $F\in {\cal
F}_1(I)$ has a decreasing iterative root $f\in {\cal F}(I)$ of odd
order $k$. Then there exists a correspondence $\kappa^{'}_f:
\{0,1,..., m\}\to \{1,...,m\}$ such that, for every $i\in \{0,1,...,
m\}$, {\rm (i)} $f^{\kappa^{'}_f(i)}(I_i)\in \Lambda(F^2)$, {\rm
(ii)} $f^{k}(I_i)\not\in \Lambda(F^2)$ for all $1\le k\le
\kappa^{'}_f(i)-1$, and {\rm (iii)} $f^{p}(I_i)\cap
f^{q}(I_i)=\emptyset$ for $1\le p\not=q\le \kappa^{'}_f(i)$.
\label{lm18}
\end{lm}

The proof of Lemma~\ref{lm18} is similar to that of
Lemma~\ref{lm14}.

By Lemma~\ref{lm18}, we can find decreasing iterative roots of $F$ on $\Lambda(F^2)$ first and then extending them to the whole domain.
The proof of Theorem~\ref{thm3} also shows that the result for square iterative roots is true for high orders. Therefore, we have

\begin{cor}\label{cor3}
Suppose that $F\in {\cal F}_1(I)$ is strictly decreasing. Then every
decreasing iterative root $f\in {\cal F}(I)$ of odd order of $F$ on
$I\backslash J(F)$ is an extension from a strictly decreasing
iterative root of the same order of $F$ on subintervals in
$\Lambda(F^2)$.
\end{cor}

Remark that strictly decreasing iterative roots of $F$ of odd order on subintervals in $\Lambda(F^2)$, where $F$ is a strictly decreasing self-mapping on $\Lambda(F^2)$, can be found in \cite[Theorem 15.8]{Kuczma68} or \cite[Theorem 11.2.4]{KCG1990}.

Similarly to the proofs of Theorems~\ref{thmD1}-\ref{thm10}, we have
the following results of decreasing iterative roots for decreasing
case.

\begin{cor}
Suppose that $F\in {\cal F}_1(I)$ is strictly decreasing and lies in
case {\bf (J1)} for a given $c_i\in J(F)$ and that $\hat{f}$ mapping
$I\setminus J(F)$ onto itself is a strictly decreasing continuous
iterative root of odd order $k$ of $F$ on $I\backslash J(F)$. Then
the multifunction
\begin{eqnarray*}\label{ddfdefineD}
f(x):=\left\{
\begin{array}{ll}
\hat{f}(x), \ \ & \ \ {\rm as}\ \ x\in I\backslash J(F),
\\
(\hat{f})^{-(k-1)}\circ F(c_i),\ \ & \ \ {\rm as}\ \  x=c_i,
\end{array}
\right.
\end{eqnarray*}
is a strictly decreasing usc
iterative root
of order $k$ of $F$ on $(I\setminus J(F))\cup \{c_i\}$.
\label{cor4}
\end{cor}

\begin{cor}
Suppose that $F\in {\cal F}_1(I)$ is strictly decreasing and lies in
case {\bf (J2)} for a given $c_i\in J(F)$ such that the set $F(c_i)$
is a compact interval and that $\hat{f}$ is a strictly decreasing
continuous iterative root of odd order $k$ of $F$ on $I\setminus
J(F)$. Then the multifunction
\begin{eqnarray*}\label{dddfdeind}
f(x):=\left\{
\begin{array}{ll}
\hat{f}(x), \ \ & \ \ {\rm as}\ \ x\in I\backslash J(F),
\\
\,[\lim_{x\to c_{i}+} \hat{f}(x),\lim_{x\to c_{i}-} \hat{f}(x)],\ \ & \ \ {\rm as}\ \  x=c_i,
\end{array}
\right.
\end{eqnarray*}
is a strictly decreasing usc
iterative root of odd order $k$ of $F$ on $(I\setminus J(F))\cup \{c_i\}$
if
$J(F)\cap (\lim_{x\to c_{i}+} \hat{f}(x),\lim_{x\to c_{i}-} \hat{f}(x))=c_i$.
\label{cor5}
\end{cor}


\section{Further discussion and remarks
}

In this section, we first give remarks on iterative roots of exclusive multifunctions.

As discussed in section 3, Theorem~\ref{thm16} is devoted to the case that the set $F(c_i)$ is a compact interval. If not, such a set $\Omega$ in the proof may not exist.
Consider $F: I:=[0,1]\to 2^{I}$ defined by
$$
F(x):=\left\{
\begin{array}{ll}
\frac{1}{16}x, \ \ & \ \ {\rm as}~~ x\in
[0,\frac{1}{2}),
\\
\,[\frac{1}{32},\frac{1}{16}],\ \ & \ \ {\rm as}~~ x = \frac{1}{2},
\\
\frac{1}{16}x+\frac{1}{32}, \ \ & \ \ {\rm as}~~ x\in (\frac{1}{2},1),
\\
\,\{\frac{3}{32},1\},\ \ & \ \ {\rm as}~~ x = 1
\end{array}
\right.
$$
for example. Clearly, $F\in {\cal F}_1(I)$. Moreover, $F$ satisfies
hypothesis ($\cal{H}$) and lies in case {\bf (J2)} with $c_i:=1$.
Note that $F(1)=\{\frac{3}{32},1\}$, is not a compact interval. For
this $F$ we cannot find a set $\Omega\subset I$ satisfying
(\ref{revisedomega}). Otherwise,
there is an $x_0\in I\backslash J(F)$ such that $\{x_0,1\}\subset \Omega$.
Further, $\Omega\subset F(1)$ implies that $x_0=\frac{3}{32}$. Then,
$
\{\frac{3}{32},1\}=F(1)=\bigcup_{j=1}^{n-1}(f_*^{j}(\frac{3}{32}))\cup\{\frac{3}{32},1\},
$
a contradiction to the fact that $f_*(\frac{3}{32})<\frac{3}{32}$
for every $f_*$ defined in (\ref{deff}).

As shown in Theorems~\ref{thm7}-\ref{thm16}, the construction of
iterative roots is simple in cases {\bf (J1)}-{\bf (J2)} because
each jump $c_i\in J(F)$ in cases {\bf (J1)}-{\bf (J2)} is
independent under iteration. In contrast, the situation of {\bf
(J3)} is much more complicated. Although Theorem~\ref{thm8}
indicates the nonexistence of strictly increasing upper
semi-continuous iterative roots of order $n>m-\ell+1$, it is still
possible for $F$ to have such a root of lower order. Let $F:
I:=[0,1]\to 2^{I}$ be defined by
$$
F(x):=\left\{
\begin{array}{ll}
\frac{1}{4}x, \ \ & \ \ {\rm as}~~ x\in
[0,\frac{1}{2}),
\\
\,[\frac{1}{8},\frac{1}{6}],\ \ & \ \ {\rm as}~~ x = \frac{1}{2},
\\
\frac{1}{6}x+\frac{1}{12}, \ \ & \ \ {\rm as}~~ x\in (\frac{1}{2},\frac{3}{4}),
\\
\,[\frac{5}{24},\frac{1}{3}],\ \ & \ \ {\rm as}~~ x = \frac{3}{4},
\\
\frac{4}{15}x+\frac{2}{15}, \ \ & \ \ {\rm as}~~ x\in (\frac{3}{4},1),
\\
\,[\frac{2}{5},\frac{1}{2}],\ \ & \ \ {\rm as}~~ x = 1.
\end{array}
\right.
$$
Obviously, $F\in {\cal F}_1(I)$. Moreover, $F$ is strictly increasing
and satisfies hypothesis ($\cal{H}$). One can check that $m=3$, and
$J(F)\cap F(1)=\{\frac{1}{2}\}$, which implies that $F$ lies in case
{\bf (J3)} with $c_i:=1$ and $\ell:=1$. By Theorem \ref{thm8}, $F$
does not have a strictly increasing usc
iterative
root of order $n>3$. However, one can verify that the mapping $f: I
\to 2^{I}$ defined by
$$
f(x):=\left\{
\begin{array}{ll}
\frac{1}{2}x, \ \ & \ \ {\rm as}~~ x\in
[0,\frac{1}{2}),
\\
\,[\frac{1}{4},\frac{1}{3}],\ \ & \ \ {\rm as}~~ x = \frac{1}{2},
\\
\frac{1}{3}x+\frac{1}{6}, \ \ & \ \ {\rm as}~~ x\in (\frac{1}{2},\frac{3}{4}),
\\
\,[\frac{5}{12},\frac{1}{2}],\ \ & \ \ {\rm as}~~ x = \frac{3}{4},
\\
\frac{4}{5}x-\frac{1}{10}, \ \ & \ \ {\rm as}~~ x\in (\frac{3}{4},1),
\\
\,[\frac{7}{10},\frac{3}{4}],\ \ & \ \ {\rm as}~~ x = 1,
\end{array}
\right.
$$
is a strictly increasing usc
iterative root of $F$ of order 2.

It is still hard to give a sufficient condition in general for those
$F$ satisfying all conditions of Theorem~\ref{thm8} to have an
increasing usc
iterative root $f$ of order
$n\leq m-\ell+1$ because the definition of $f$ at the specific jump
$c_i$ depends on the values of $f$ at other jumps such that
$J(F)\cap f^n(c_i)=\{c_{k_1}, c_{k_2}, \ldots, c_{k_\ell}\}$.
Actually, for $n=2$, the simplest case, an iterative root $f$ needs
to satisfy $J(F)\cap f^2(c_i)=\{c_{k_1}, c_{k_2}, \ldots,
c_{k_\ell}\}$, which requires that $\{c_{k_1}, c_{k_2}, \ldots,
c_{k_\ell}\}\subset f(f(c_i))$. The difficulty is to find suitable
jumps $c^*\in J(F)$ such that $c^*\in f(c_i)$ and $\{c_{k_1},
c_{k_2}, \ldots, c_{k_\ell}\}\subset f(c^*)$. Hence, the definition
of $f$ at the specific $c_i$, depends on the values of $f$ at those
jumps $c^*$. However, we can give some properties for those roots:
If $f\in {\cal F}(I)$ is an increasing iterative root of $F$ of
order $n\leq m-\ell+1$, where $F$ satisfies conditions of
Theorem~\ref{thm8}, then
\begin{eqnarray}
f(c_i)\cap J(F)= S_1 ,~ f(S_j)\cap J(F)= S_{j+1}\cap f(J(F)),~j=1,...,n-1,
\label{inclusions}
\end{eqnarray}
and there is an integer $\tau\in \{1,...,m\}$ such that
\begin{eqnarray}
c_{i} \xrightarrow[]{f} N_{i, {1}} \xrightarrow[]{f} N_{i, {2}}
\xrightarrow[]{f}... \xrightarrow[]{f} N_{i, {\tau}}
\xrightarrow[]{f}  K,
\label{fci}
\end{eqnarray}
where
$N_{i, k}:=f^k(c_i)\setminus J(F)$, $k=1, 2,\ldots,\tau$.
In fact,
$
f(f^j(c_i)\cap J(F))\cap J(F)=S_{j+1}\cap f(J(F))
$
for $j\in \{1,2,\ldots, n-1\}$, which implies (\ref{inclusions}) by (\ref{sf}) and (\ref{fci}) is proved by Lemma~\ref{lm16}.

The problem of iterative roots in case {\bf (J4)} is still open.
The main reason again is that complicated computations are caused by
the definition of $f$ at the specific jump $c_i$, which depends on
the values of $f$ at other jumps such that $J(F)\cap f^n(c_i)=\{c_i,c_{k_1}, c_{k_2}, \ldots, c_{k_\ell}\}$,
the same as for iterative roots of lower order in case {\bf (J3)}. Notice that the difference between case {\bf (J4)} and case {\bf (J3)}
is that $c_i\in f(c_i)$ in case {\bf (J3)} ($c_i\not\in f(c_i)$ if $F$ lies in case {\bf (J4)}).
However, the existence of iterative roots in case {\bf (J4)} is possible.
Let $F: I:=[0,1]\to 2^{I}$ be defined by
$$
F(x):=\left\{
\begin{array}{ll}
\frac{1}{4}x, \ \ & \ \ {\rm as}~~ x\in
[0,\frac{1}{2}),
\\
\,[\frac{1}{8},\frac{1}{6}],\ \ & \ \ {\rm as}~~ x = \frac{1}{2},
\\
\frac{1}{6}x+\frac{1}{12}, \ \ & \ \ {\rm as}~~ x\in (\frac{1}{2},1),
\\
\,[\frac{1}{4},1],\ \ & \ \ {\rm as}~~ x = 1.
\end{array}
\right.
$$
Obviously, $F\in {\cal F}_1(I)$ is strictly increasing and lies in
case {\bf (J4)} with $c_i:=1$ and $\ell:=1$. One can check that the
mapping $f: I \to 2^{I}$ defined by
$$
f(x):=\left\{
\begin{array}{ll}
\frac{1}{2}x, \ \ & \ \ {\rm as}~~ x\in
[0,\frac{1}{2}),
\\
\,[\frac{1}{4},\frac{1}{3}],\ \ & \ \ {\rm as}~~ x = \frac{1}{2},
\\
\frac{1}{3}x+\frac{1}{6}, \ \ & \ \ {\rm as}~~ x\in (\frac{1}{2},1),
\\
\,[\frac{1}{2},1],\ \ & \ \ {\rm as}~~ x = 1,
\end{array}
\right.
$$
is a strictly increasing usc
iterative root of $F$ of order 2.

Compared with results in reference \cite{Li-Jar-Zhang09Publ},
increasing iterative roots of order 2 were found for those strictly
increasing functions $F\in {\cal F}_1(I)$, each of which has a unique
jump $c\in I$. Theorem 5 of \cite{Li-Jar-Zhang09Publ} requires that
$F(x)\neq \{c\}$ for all $x\in I\backslash \{c\}$ and $c\not\in
F(c)$, which implies that $\zeta(F)=1$, i.e., $F\in {\cal
F}_1(I)$, and that $F$ lies in our case {\bf (J1)}. Since our
Theorems~\ref{thm7}-\ref{thm8} generally consider roots of order $n\geq 2$
for multifunctions having more than one jumps, the cases {\bf (J3)}-{\bf (J4)} are not encountered in \cite{Li-Jar-Zhang09Publ}. Moreover, the
case {\bf (J2)}, which is dealt with in Theorem~\ref{thm16}, was
considered in Theorem 6 of \cite{Li-Jar-Zhang09Publ} with a set $M$.
Our Theorem~\ref{thm16} gives $f(x)=\Omega:=[\lim_{s\to
c_{i}-}f_*(s),c_i]$ for $x=c_i$, which is actually a formulation of
the set $M$.

The problem of strictly decreasing iterative roots of a general
order $n$ for exclusive multifunctions is reduced to finding strictly decreasing iterative roots
of order $2$, as showed in section 4. In case {\bf (J1)}, unlike Theorem 5 of
\cite{Li-Jar-Zhang09Publ}, where all decreasing semi-continuous
square iterative roots
are found for those strictly increasing functions $F\in {\cal F}_1(I)$
each of which has a unique jump, our Theorem~\ref{thmD1} considers
multifunctions having more than one jumps. In case {\bf (J2)}, our
Theorem~\ref{thm10} requires (\ref{Jf^*}), i.e., $J(F)\cap
(\lim_{x\to c_{i}+} f^*(x),\lim_{x\to c_{i}-} f^*(x))=c_i$, which
replaces the hypothesis ($\cal{H}$) used in Theorem~\ref{thm16},
because Lemma~\ref{lm17} does not hold for decreasing multifunction
anymore. Our Theorem~\ref{thm10} gives $f(x)=[\lim_{s\to c_{i}+}
f^*(s),\lim_{s\to c_{i}-} f^*(s)]$ for $x=c_i$, which is actually a
formulation of the set $M$ given in Theorem 6 of
\cite{Li-Jar-Zhang09Publ}.

From Theorems~\ref{thmD1}-\ref{thm10}, we can construct
decreasing iterative roots in cases {\bf (J1)} -{\bf (J2)} because
every jump $c_i\in J(F)$ in those cases is independent under
iteration. In contrast, in cases {\bf (J3)} -{\bf (J4)} the problem
of decreasing iterative roots remains difficult because the
decreasing monotonicity reverses the orientation and the value of a
decreasing iterative root $f$ at the specific jump $c_i$ depends on
the values of $f$ at other jumps.

As above we focus on those multifunctions having intensity 1, as
indicated at the end of section 2. However, there are many
multifunctions with intensity larger than 1. Although the case of
$\zeta(F)>1$ is more complicated, similarly to Theorem 1 of
\cite{ZWN-Pol} for single-valued functions, we have a basic result:
{\it Any $F\in {\cal F}(I)$ with intensity $\zeta(F)>1$ has no upper
semi-continuous iterative roots of order $n>\#J(F)$. }For an
indirect proof, assume that $F$ has a
usc
iterative root $f$ of order $n>\#J(F)$. Since $\zeta(F)>1$, i.e.,
$\#J(F)<\#J(F^2)$, we have
$$
0<\#J(f)<\#J(f^2)<...<\#J(f^n)=\#J(F)<\#J(F^2)
$$
by Lemma~\ref{lm7}. It follows that $\#J(F)\ge n$, a contradiction
to the fact that $\#J(F)<n$.
A corollary of the basic result is:
{\it
Any $F\in {\cal F}(I)$ with a unique jump and $\zeta(F)>1$ has no usc
iterative roots of any order $n\ge 2$.
}


In this paper we mainly consider strictly monotone usc multifunctions. As known in the
theory (e.g. \cite{Li-Yang-Zhang08,LiuJarczykLiZhang2011NA,Liu-Zhang11}) of
iterative roots for single-valued functions, the case without
monotonicity is much more complicated. So it is more challenging to
discuss on iterative roots for those usc
multifunctions without monotonicity.
\\
\\
{\bf Acknowledgements}
\\
\\
The authors are very grateful to the reviewers for their carefully checking and helpful suggestions.


\end{document}